# BIBLIOGRAPHY

## Publications of J. L. Doob

### 1932

[1] The boundary values of analytic functions. *Trans. Amer. Math. Soc.* **34** 153–170.

[2] On a theorem of Gross and Iverson. *Ann. of Math.* **33** 753–757.

### 1933

[3] Systems of algebraic difference equations (with J. F. Ritt). *Amer. J. Math.* **55** 505–514.

### 1934

[4] The boundary values of analytic functions II. *Trans. Amer. Math. Soc.* **35** 418–451.

[5] Stochastic processes and statistics. *Proc. Natl. Acad. Sci. USA* **20** 376–379.

[6] On analytic functions with positive imaginary parts (with B. O. Koopman). *Bull. Amer. Math. Soc.* **40** 601–605.

[7] Probability and statistics. *Trans. Amer. Math. Soc.* **36** 759–775.

### 1935

[8] The limiting distributions of certain statistics. *Ann. Math. Statist.* **6** 160–169.

[9] The ranges of analytic functions. *Ann. of Math.* **36** 117–126.

### 1936

[10] Note on probability. *Ann. of Math.* **37** 363–367.

[11] Statistical estimation. *Trans. Amer. Math. Soc.* **39** 410–421.

---